\documentclass[11pt]{article}
\usepackage{amsmath,amssymb}

\def\emptyset{\varnothing}

\def\kappa{\varkappa}

\let\epsilon\varepsilon

\let\phi\varphi

\def\n{\noindent}
\def\upn{}
\def\blanksquare{\,\,\,$\sqcup\!\!\!\!\sqcap$}         
\def\qed{\hfill\blanksquare\linebreak\smallskip\par}   
\def\bline(#1,#2)(#3,#4)(#5){\put(#1,#2){\line(#3,#4){#5}}} 

\begin{document}


\title{On Raw Coding of Chaotic Dynamics}

\author{{\Large Michael Blank}\thanks{Institute for Information 
        Transmission Problems, Moscow  {\tt blank@iitp.ru}} 
        \thanks{Supported in part by the Russian Foundation for 
        Basic Research, project no.~05-01-00449, 
        Grant of the President of the Russian Federation for 
        Leading Scientific Schools, no.~934.2003.1, 
        and the French Ministry of Education.}}

\date{March 24, 2006} 

\maketitle

\begin{abstract}
We study raw coding of trajectories of a chaotic dynamical system by
sequences of elements from a finite alphabet and show that there is a
fundamental constraint on differences between codes corresponding to
different trajectories of the dynamical system.
\end{abstract}

\section{Introduction}\label{s:intro}

By a \emph{raw coding\/} of a chaotic discrete time dynamical system
$(T,X,\mathcal{B},\mu)$ (defined by a measurable map $T$ of a measurable
space $(X,\mathcal{B},\mu)$ into itself with a $\sigma$-algebra of measurable
sets $\mathcal{B}$ and a $T$\nobreakdash-invariant probability measure
$\mu$), we mean a representation of trajectories
$\{T^tx\}\strut_{t\in\mathbb{Z}_+}$ of this system as sequences of elements
from a finite alphabet $\mathcal{A}$. In other words, one defines a mapping
$\Xi:X\to\mathcal{A}$, which takes points of the original phase space $X$ of
the dynamical system to elements of the alphabet $\mathcal{A}$. This mapping
defines a partition $\xi:=\{X_1,X_2,\dots,X_M\}$ of the phase space $X$ into
disjoint measurable subsets $X_i\in\mathcal{B}$, $i\in\{1,2,\dots,M\}$, with
a subsequent encoding of a trajectory by a sequence of numbers corresponding
to elements of the partition which the trajectory successively visits; i.e.,
$\Xi(x):=i$ if $x\in X_i$. In some cases (if there exists a finite Markov
partition; see, e.g,. \cite{Bo,Sinai}), information about raw encoded
trajectories allows one to completely reconstruct all topological
characteristics of the system under consideration. Note however that, even if
the existence of a finite Markov partition is proved rigorously, its
practicable application is rather questionable due to the high instability of
the chaotic map $T$. Moreover, it might happen that an approximate Markov
partition has even worse statistical characteristics (i.e., more differing
from the characteristics of the original system) than an ``arbitrary'' one
\cite{BK}.

The aim of the present paper is to demonstrate that, under a rather weak
assumption about a deterministic dynamical system (weak mixing condition),
there is a fundamental constraint on elementwise differences between codes
corresponding to different trajectories. Using standard definitions and
constructions of ergodic theory, the main of which will be described in
Section~2 and whose detailed analysis can be found, e.g., in \cite{Sinai,
Bl}, this statement can be formulated as follows.

\medskip\n{\bf Theorem.} \label{t:variability}
Let a dynamical system\/ $(T,X,\mathcal{B},\mu)$ satisfy the weak mixing
property. Consider a finite measurable partition\/
$\xi:=\{X_1,X_2,\dots,X_M\}$ with\/ $\prod\limits_{i=1}^M\mu(X_i)>0$.
Then\upn, for any positive integers $N,L<\infty$ and for almost every\/
\upn(with respect to the direct product measure $\mu^N$\upn) collection
$\bar{x}:=\{x_1,x_2,\dots,x_N\}\in X^N$\upn, there exists a time moment\/
$t_0\ge0$ such that\/ $\Xi(T^{t_0+t}x_i)=\Xi(T^{t_0+t}x_j)$ for~all\/
$t\in\{1,\dots,L\}$\upn, $i,j\in\{1,\dots,N\}$.
\bigskip

In other words, all of these $N$ codes of different trajectories
simultaneously contain arbitrarily long subsequences coinciding both in space
and time. In Section~3, after the proof of this result, we discuss possible
weakenings of the assumptions under which it holds.

From the point of view of numerous applications, let us mention connections
of this problem to the analysis of properties of random number generators.
Let $\{x_i^t\}\strut_{t\in\mathbb{Z}_+}$, $i\in\{1,\dots,N\}$,
be~$N$~different sequences of pseudorandom numbers obtained from the same
random number generator. We assume that those pseudorandom numbers belong to
a finite alphabet $\mathcal{A}$. It is desirable that, despite a proximity
(ideally, coincidence) of statistical properties of these realizations, their
pointwise differences should be as large as possible. The theorem shows that
the pointwise differences between different realizations are fundamentally
constrained. Although this claim is somewhat unexpected, it completely agrees
with well-known properties of Bernoulli or, more generally, Markov random
sequences. Indeed, consider a Bernoulli random sequence over an alphabet with
two elements $\{0,1\}$, taken with probabilities $p$ and $q=1-p$,
respectively. Consider an arbitrary finite sequence $\bar{a}$ over this
alphabet, with $K\le L$ elements $0$ and $L-K$ elements $1$. For any given
positive integer $t_0$, the probability that a realization of this Bernoulli
random sequence coincides with $\bar{a}$ on the time interval from $t_0$ to
$t_0+L-1$ is $p^Kq^{L-K}>0$; it does not depend on an initial moment $t_0$.
Thus, for an arbitrary finite number of realizations of this Bernoulli
sequence, the claim of the theorem holds true with probability one. One can
similarly prove this statement for a Markov random chain with a finite number
of states $A=\{1,\dots,M\}$ and a transition probability matrix $\pi$
satisfying the assumption $\pi^\kappa>0$ for some $\kappa\in\mathbb{Z}_+$
(this is a direct analog of the weak mixing property for a Markov chain with
a finite number of states).

It is worth to note connections of this problem to the Ulam finite Markov
chain approximation scheme and to the so-called Bowen specification property
(see, e.g. \cite{Bl,Bo}) for chaotic dynamical systems.

\section{Necessary definitions and constructions}\label{s:defs}

Here we give a short description of standard definitions and constructions
from ergodic theory which we need to prove the theorem.

Recall that a measure $\mu$ is \emph{$T$-invariant\/} if and only if
$$
\int\phi\circ T\,d\mu=\int\phi\,d\mu
$$
for any $\mu$-integrable function $\phi\colon X\to\mathbb{R}^1$.

A measurable function $\phi\colon X\to\mathbb{R}^1$ is called
\emph{invariant\/} with respect to a dynamical system $(T,X,\mathcal{B},\mu)$
(or simply $T$-invariant) if $\phi=\phi\circ T$ almost everywhere with
respect to the measure $\mu$.

A dynamical system $(T,X,\mathcal{B},\mu)$ is \emph{ergodic\/} if each
$T$-invariant function is a constant $\mu$-a.e.

A dynamical system $(T,X,\mathcal{B},\mu)$ is \emph{weakly mixing\/} if
$$
\frac1n\sum_{k=0}^{n-1}\bigl|\mu(T^{-k}A\cap B) - \mu(A)\mu(B)\bigr|
\xrightarrow{n\to\infty\,}0, \quad \forall A,B\in\mathcal{B}.
$$

A \emph{direct product\/} of two dynamical systems
$(T',X',\mathcal{B}',\mu')$ and $(T'',X'',\mathcal{B}'',\mu'')$ is a new
dynamical system $(T'\otimes
T'',X'\otimes{X''},\mathcal{B}'\otimes\mathcal{B}'',\mu'\otimes\mu'')$, where
the map $T'\otimes T''\colon X'\otimes{X''}\to X'\otimes{X''}$ is defined by
the relation
$$
T'\otimes T''(x',x''):=(T'x',T''x''),
$$
while all other objects are standard direct product of spaces,
$\sigma$-algebras, and measures, respectively.

By $A^N$ we denote the $N$th power of the set $A$, i.e., the direct product
of $N$ identical sets $A\in\mathcal{B}$, and by
$(T^{\otimes{N}},X^N,\mathcal{B}^N,\mu^N)$ the direct product of $N$
identical dynamical systems $(T,X,\mathcal{B},\mu)$.

A collection $\xi:=\{X_1,X_2,\dots,X_M\}$ is called a \emph{measurable
partition\/} of a measurable space $(X,\mathcal{B},\mu)$ if
$X_i\in\mathcal{B} ~\forall i$, $X_i\cap X_j=\emptyset ~\forall i\ne j$, and
$\bigcup\limits_i X_i=X$.

\section{Proof of the theorem}\label{s:proof}

The proof of this result consists of two steps. First we show that a
sufficient condition for the claim of the theorem is the ergodicity of the
direct product of $N$ copies of the original dynamical system. Then we
demonstrate that the latter property follows from the conditions of the
theorem.

Denote $\xi^{(0)}:=\xi$ and inductively define a sequence of measurable
partitions $\xi^{(n)}$ with $n\in\mathbb{Z}_+$ by the following relation:
$$
\xi^{(n+1)} := \xi^{(n)} \cap T^{-1}\xi^{(n)}.
$$
It is immediate to check that the constructed collections of sets $\xi^{(n)}$
are measurable partitions of the space $(X,\mathcal{B},\mu)$ for any positive
integer $n$. The partition $\xi^{(n)}$ is called the $n$th
\emph{refinement\/} of the partition $\xi$.

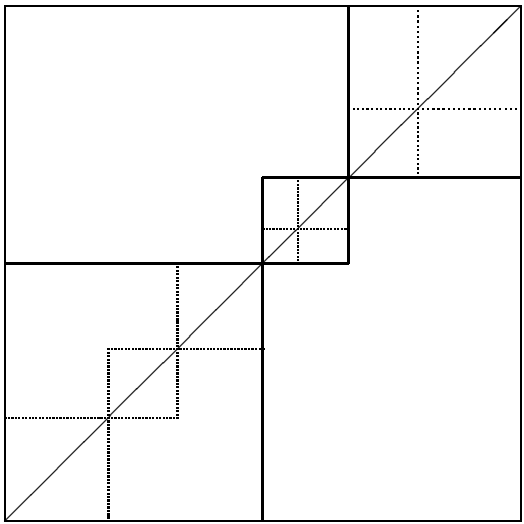
\begin{figure}[tp]
\centering\unitlength=1.3pt
\begin{picture}(150,150)
\bline(0,0)(1,0)(150) \bline(0,0)(0,1)(150) \bline(0,150)(1,0)(150)
\bline(150,0)(0,1)(150) \bline(0,0)(1,1)(150) \thicklines
\bline(0,75)(1,0)(75) \bline(75,0)(0,1)(75) \bline(75,75)(1,0)(25)
\bline(75,75)(0,1)(25) \bline(100,100)(-1,0)(25) \bline(100,100)(0,-1)(25)
\bline(100,100)(1,0)(50) \bline(100,100)(0,1)(50) \thinlines
\bezier{30}(30,30)(30,15)(30,0) \bezier{30}(30,30)(15,30)(0,30)
\bezier{20}(30,30)(30,40)(30,50) \bezier{20}(30,30)(40,30)(50,30)
\bezier{20}(50,50)(50,40)(50,30) \bezier{20}(50,50)(40,50)(30,50)
\bezier{25}(50,50)(50,62)(50,75) \bezier{25}(50,50)(62,50)(75,50)
\bezier{10}(85,85)(85,80)(85,75) \bezier{10}(85,85)(80,85)(75,85)
\bezier{15}(85,85)(85,92)(85,100) \bezier{15}(85,85)(92,85)(100,85)
\bezier{15}(120,120)(120,110)(120,100) \bezier{15}(120,120)(110,120)(100,120)
\bezier{20}(120,120)(120,135)(120,150) \bezier{20}(120,120)(135,120)(150,120)
\end{picture}
\caption{``Toothed neighborhoods of the diagonal'' $[0,1]^2$ created by
collections of the sets $(\xi^{(n)})^2$ (the~boundary is indicated by solid
lines) and $(\xi^{(n+1)})^2$ (the boundary is indicated by dashed
lines).}\label{f:toothed-neighborhood}
\end{figure}

By means of these partitions, for each $n\in\mathbb{Z}_+$ we define a
collection of sets
$$
(\xi^{(n)})^N:=\bigcup_{\eta\in\xi^{(n)}}\eta^N,
$$
which are ``toothed neighborhoods of the diagonal'' in the direct product
$X^N$ (see Fig.~\ref{f:toothed-neighborhood}).

Observe that the inclusion $\bar{x}:=\{x_1,\dots,x_N\}\in\eta^N$ with
$\eta\in\xi^{(n)}$ implies the inclusion
$$
(T^{\otimes N})^t\bar{x} \in (T^{\otimes N})^t\eta^N \in(\xi^{n-t})^N
$$
for all $t\in\{0,1,\dots,n\}$; hence,
$$
\Xi(T^tx_i) = \Xi(T^tx_j), \quad \forall i,j\in\{1,2,\dots,N\},\quad
t\in\{0,1,\dots,n\}.
$$

Therefore, if we prove that, for $\mu^N$-a.a.\ collection of points
$\bar{x}:=\{x_1,\dots,x_N\}\in X^N$, there exists a time moment $t_0\ge0$
such that
$$
(T^{\otimes N})^{t_0}\bar{x} \in\eta^N\in(\xi^{(L)})^N,
$$
then the claim of the theorem will follow.

Let a dynamical system $(\tau,Y,\mathcal{B}_Y,\nu)$ be ergodic. Then, for any
pair of measurable sets $A,B\in\mathcal{B}_Y$ with $\nu(A)\nu(B)>0$, there
exists a positive integer $\kappa=\kappa(A,B)<\infty$ such that $\tau^\kappa
A\cap B\ne\emptyset$. Indeed, assume that this is not true, i.e.,
$\tau^nA\cap B=\emptyset$ for any positive integer $n$. Consider a measurable
set
$$
A_\infty:=\bigcup_{n\in\mathbb{Z}_+}\tau^nA.
$$
Then $\nu(A_\infty)\ge\nu(A)>0$ and $A_\infty\cap B=\emptyset$. Therefore,
the indicator function of a measurable set~$A_\infty$ of positive
$\nu$-measure is $\tau$-invariant but is not a constant a.e., which
contradicts the ergodicity.

Thus, it suffices to show that the dynamical system $(T^{\otimes
N},X^N,\mathcal{B}^N,\mu^N)$ is ergodic. For that, we use the fact that the
weak mixing property is preserved under the direct product of weakly mixing
dynamical systems (see e.g. \cite{Sinai}). To complete the proof, its remains
to note that weak mixing implies ergodicity. \qed

\smallskip
It is interesting that one cannot weaken the conditions of the theorem by
changing the weak mixing of the original dynamical system to the ergodicity.
The problem is that the direct product of ergodic dynamical systems need not
also be ergodic: consider a direct product of two identical irrational unit
circle rotations.
\begin{figure}[tp]
\centering\unitlength=1.3pt
\begin{picture}(150,150)
\bline(0,0)(1,0)(150) \bline(0,0)(0,1)(150) \bline(0,150)(1,0)(150)
\bline(150,0)(0,1)(150) \bline(0,0)(1,1)(150)
\bezier{50}(0,75)(75,75)(150,75) \thicklines \put(0,0){\vector(1,2){37}}
\put(37,0){\vector(1,2){75}} \bline(113,75)(1,2)(37) \put(-5,-10){$0$}
\put(150,-10){$1$} \put(-5,145){$1$} \put(-16,35){$X_L$} \put(-16,115){$X_R$}
\end{picture}
\vskip13pt \caption{Counterexample to the necessity of the weak mixing
property.}\label{f:counter}
\end{figure}
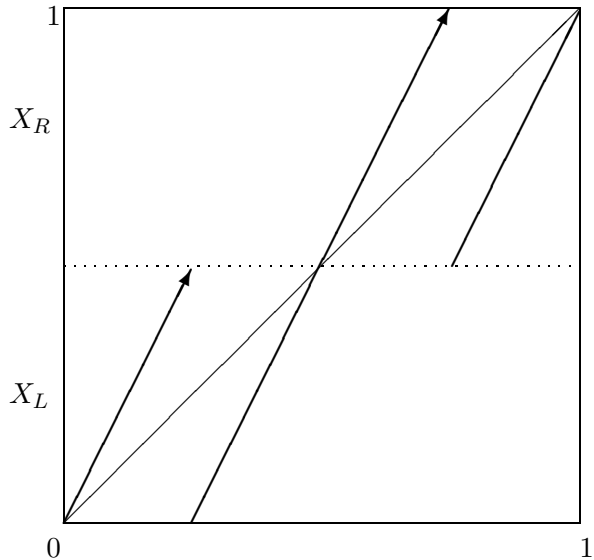

On the other hand, the weak mixing condition is not necessary either. We
shall show that, for some nonergodic dynamical systems, there exist
partitions for which the claim of the theorem holds true. Indeed, let $T$ be
a mapping from the unit interval $X=[0,1]$ into itself defined by the
relation
\begin{equation}\label{e:nonerg}
T x :=
\begin{cases}
2x & \text{if}\ 0\le x <1/4,\\ 2x-1/2 & \text{if}\ 1/4\le x <3/4,\\ 2x-1 &
\text{otherwise}.
\end{cases}
\end{equation}
The graph of this mapping is shown in Fig.~\ref{f:counter}. It is immediate
to check that the restriction of the map $T$ to each of the half-intervals
$X_L:=[0,1/2]$ (left) and $X_R:=[1/2,1]$ (right) with the normalized Lebesgue
measure $m$ is a weakly (and even strongly) mixing dynamical system. Thus,
the Lebesgue invariant measure of the map $T$ considered on the entire unit
interval is not ergodic. Consider a partition $\xi$ with
$\prod\limits_{i=1}^M\mu(X_i)>0$ satisfying the assumption that
$X_1:=[1/2-2^{-k},1/2+2^{-k})$ for some $k\in\mathbb{Z}_+$. Denote by $\xi_L$
and $\xi_R$ the restrictions of the partition $\xi$ to the left and right
half-intervals, respectively. The claim of the theorem is satisfied for each
of the dynamical systems $(T,X_L,\mathcal{B},m)$ and $(T,X_R,\mathcal{B},m)$
equipped with the partitions $\xi_L$ and $\xi_R$, respectively. Moreover,
it~turns out that, for the nonergodic dynamical system $(T,X,\mathcal{B},m)$
with the partition $\xi$, the claim of the theorem holds as well.

To check this, it is sufficient to show that, for any measurable set
$A\subseteq[0,1]^N$ of positive Lebesgue measure, there exists an integer
$t_0$ such that the set $(T^{\otimes N})^{t_0}A$ has an intersection of
positive Lebesgue measure with the ``toothed neighborhood of the diagonal''
in $[0,1]^N$ generated by the $N$th power of the $L$th refinement of the
partition $\xi$. Observe that, by the construction, each of these refinements
has an element containing an open neighborhood of the point $1/2$. Denote by
$B\in[0,1]^N$ the $N$th power of this open neighborhood contained in an
element of the partition~$\xi^{(L)}$. We need to show that
\begin{equation}\label{e:intersec}
(T^{\otimes N})^{t_0}A\cap B\ne\emptyset
\end{equation}
for some positive integer $t_0$.

By a \emph{quadrant}, we shall call a direct product of $n\in\{0,1,\dots,N\}$
copies of the interval $X_L$ and $N-n$ copies of the interval $X_R$ taken in
an arbitrary order. Clearly, the union of all possible quadrants coincides
with $[0,1]^N$. On the other hand, each of the quadrants is invariant under
the map $T^{\otimes N}$. Moreover, the restriction of the dynamical system
$(T^{\otimes N},X^N,\mathcal{B}^N,\mu^N)$ to each of the quadrants is again a
direct product of mixing dynamical systems. Now there exists a quadrant whose
intersection with the set $A$ is of positive Lebesgue measure. Denote this
quadrant by $\Delta$ and~set $A_\Delta:=A\cap\Delta$ and
$B_\Delta:=B\cap\Delta$. Observe that, by the construction, we have
$m(A_\Delta)m(B_\Delta)>0$. Therefore, using the same argument as in the
proof of the theorem, we get
$$
(T^{\otimes N})^{t_0}A_\Delta\cap B_\Delta\ne\emptyset
$$
for some $t_0>0$, which implies (\ref{e:intersec}).

Therefore, the claim of the theorem about the existence of arbitrarily long
coinciding code segments takes place for this nonergodic dynamical system. In
a sense, the element $X_1$ of the partition~$\xi$ plays a role of a
``bridge'' between ergodic components of the nonergodic dynamical system
$(T,[0,1],\mathcal{B},m)$.

\medskip
The author thanks D.~Buzzi for useful discussions.

\end{document}